\documentclass[12pt,fullpage,leqno]{article}

\usepackage{amsfonts}
\usepackage{amsmath}
\usepackage{epsfig}
\usepackage{graphicx}
\usepackage{xcolor}

\def\R{\hbox{\bf R}}

\def\T{{\mathbb T}}

\newcommand{\ba}{\begin{eqnarray}}
\newcommand{\ea}{\end{eqnarray}}

\newtheorem{theo}{\bf Theorem}[section]
\newtheorem{lem}[theo]{\bf Lemma}

\newtheorem{defi}[theo]{\bf Definition}

\renewcommand{\R}{{\mathbb R}}

\newenvironment{Proofc}[1]{\smallskip\par\noindent\textsc{#1}\quad}%
  {\hfill$\Box$\bigskip\par}

\nonstopmode

\begin{document}

\title{\bf A new contraction family\\ for porous medium\\ and fast diffusion equations}
\author{G. Chmaycem$^{b,*}$, M. Jazar\footnote{LaMA-Liban, Lebanese University, P.O. Box 37 Tripoli,
Lebanon. E-mail: ghada.chmaycem@gmail.com (G. Chmaycem),
mjazar@laser-lb.org (M. Jazar).
\newline \indent $\,\,{}^{b}$ Universit\'{e} Paris-Est,
CERMICS, Ecole des Ponts ParisTech, 6 et 8 avenue Blaise Pascal,
Cit\'e Descartes Champs-sur-Marne, 77455 Marne-la-Vall\'ee Cedex 2,
France. E-mail: chemaycg@cermics.enpc.fr (G. Chmaycem),
monneau@cermics.enpc.fr (R. Monneau)} , R. Monneau$^{b}$}
\date{June 14, 2014}
\maketitle


\centerline{\small{\bf{Abstract}}}
{\small{In this paper, we present a surprising  two-dimensional contraction family for porous medium and fast diffusion equations.
This approach provides new a priori estimates on the solutions, even for the standard heat equation.
}}\hfill\break

 \noindent{\small{\bf{AMS Classification:}}} {\small{{35K55, 35K65}.}}\hfill\break
 \noindent{\small{\bf{Keywords:}}} {\small{contraction, porous medium equation, fast diffusion.
 }}\hfill\break


\section{Introduction}

In this paper, we answer a long standing open question about the existence of new contractions for porous medium type equations.
For $m>0$ and $d\ge 1$, we consider 
nonnegative solutions $U(t,x)$ of  the following normalized equation
\begin{equation}\label{eq::1}
m U_t = \Delta U^m \quad \mbox{on}\quad (0,+\infty)\times \Omega= Q \quad \mbox{with}\quad \Omega=\R^d.
\end{equation}
For $m>1$, this equation is called the porous medium equation (PME), while for $m<1$, it is called the fast diffusion equation (FDE).
We choose the following standard definition of solution:
\begin{defi}\label{defi::1}{\bf (Notion of solution)}\\
Let $m>0$.  We say that $U\ge 0$ is a solution of (\ref{eq::1}) with intial data $0\le U_0\in L^1(\Omega)$, 
if $U, U^m\in L^1_{loc}(Q)$ and for every function $\eta\in C_c^\infty(\overline{Q})$ with compact support, we have
$$\int\hspace{-0,3cm}\int_Q \left\{\frac{U^m}{m}\Delta \eta + U\eta_t\right\}\ dxdt+ \int_\Omega U_0(x)\eta(0,x) \ dx= 0.$$
\end{defi}
We recall that for $m>0$ and any initial data satisfying
\begin{equation}\label{eq::r1}
U_0\in L^1(\Omega)\cap L^{\infty}(\Omega)\quad \mbox{with}\quad U_0\ge 0,
\end{equation}
there is exactly one solution $U\in C([0,+\infty); L^1(\Omega))\cap L^{\infty}(Q)$ in the sense of Definition \ref{defi::1} (see \cite{V1,V2}).
In the special case where $m>m_c(d):=\max\left(0, \frac{d-2}{d}\right)$, it is well known that it is possible to remove the condition $U_0\in L^\infty(\Omega)$ to get existence and uniqueness in some weaker space. To present our results in a concise and  unified way for $m>0$, we will stick under assumptions (\ref{eq::r1}).

In any dimensions, it is known (see \cite{V} and the references therein), that this equation has at least three types of contractions:
in $L^1$, in $H^{-1}$ and finally for the $2$-Wasserstein distance if $m>m_c(d)$. This last contraction property has been discovered by Otto in \cite{O2} and published in \cite{O1} (see also later \cite{CMV}). 
It is also known for the PME, that for the $p$-Wasserstein distance, there is no contraction for $d\ge 2$ and $p>p_1(m,d)$, while there is contraction for any $p\in [1,+\infty]$ in dimension $d=1$ (see \cite{V}).
Note also that in \cite{V2} (Theorem A.5 p. 583), the author proves that the PME is not contractive with respect to $L^p$ norm for $m\ge 2$ if $p>p_2(m,d)$.

\begin{figure}[htb]
\begin{center}
\includegraphics[width=0.8\textwidth]{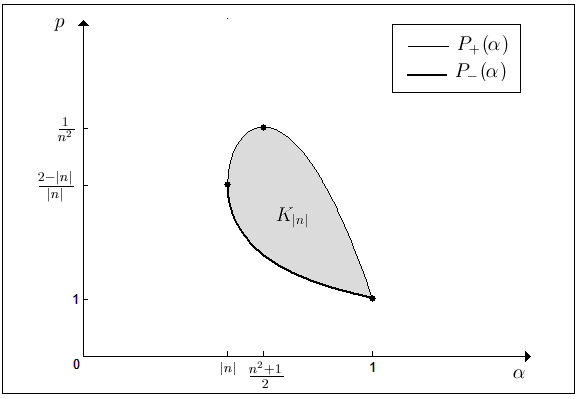}
\caption{The set $K_{|n|}$. }\label{fig-Kn}
\end{center}
\end{figure}

In this paper, for $0<m<2$, we present a new family of contractions for this equation in any dimensions, which extends the $L^1$ contraction properties. Our contraction can be seen as the fourth known contraction for this equation.
Even for the case $m=1$, our approach leads to new results for the standard heat equation.\\
More precisely, for $U$ and $V$ two nonnegative solutions of (\ref{eq::1}), we show that the following quantity
\begin{equation}\label{eq::r2}
\int_{\Omega}|U^\alpha-V^\alpha|^p
\end{equation}
is a Lyapunov functional which is nonincreasing in time for all $(\alpha,p)$ in some admissible set.\\
For convenience, we will work in the whole paper with
$$n=m-1.$$
For $0<|n|<1$, we define the admissible (convex) set (which is skeched on Figure \ref{fig-Kn})

\begin{equation}\label{eq::2}
K_{|n|}=\left\{(\alpha,p)\in\left[|n|,1\right]\times {\mathbb{R}}; \quad 
P_-(\alpha)\leq p\leq P_+(\alpha),\quad \mbox{with}\quad P_\pm(\alpha)=1+\frac{2}{n^2}(1-\alpha)(\alpha\pm \sqrt{\alpha^2-n^2})\right\}.
\end{equation}
For $n=0$, we also set $K_0=\left\{(\alpha,p)\in\left(0,1\right]\times {\mathbb{R}}; \quad  \alpha p \ge 1\right\}$.
One may wonder if contraction (\ref{eq::r2}) is related or not to some gradient flow structure of the equation.
Indeed, for a given positive solution $U$ of (\ref{eq::1}), we can set $u=U^\alpha$
which solves
\begin{equation}\label{eq::1-bis}
u_t = u^{\gamma-1}\left( u\Delta u +\bar \gamma\vert \nabla u\vert^2\right) \quad \mbox{with} \quad \gamma = \frac{n}{\alpha}\end{equation}
and $\bar \gamma = \gamma-1+\alpha^{-1}$.
In the special case where   $\bar \gamma =\gamma/2$, 
it can be seen that equation (\ref{eq::1-bis}) is the negative $L^2$-gradient flow of some energy, i.e. it solves
$$u_t = -\nabla_{L^2} E(u) \quad\quad  \mbox{with}\quad \quad E(u)=\int_\Omega u^\gamma \frac{|\nabla u|^2}{2}.$$
Moreover, it is easy to check (computing the hessian of $E$) that $E$ is convex if $\gamma\in [-1,0]$. This corresponds exactly to the points $(\alpha,p)=(1+\frac{n}{2},2)\in K_{|n|}$ that we capture for $n\in [-\frac23,0]$.
This interpretation is similar to the derivation of Yamabe flow (see for instance \cite{Y}) with the difference here that the exponent of the fast diffusion is not directly related to the space dimension.
Nevertheless, except this very exceptional case, contraction (\ref{eq::r2}) does not seem to be related to any gradient flow structure.

The paper is organized as follows. 
In Section 2, we present our main results. In Section 3, we prove the new contraction property.
Finally in Section 4, we prove several gradient decay estimates.

\section{Main results}

We state precisely our main results.
With the previous notation, and for $u,v>0$, we define the symmetric matrix for $p>1$ and $\alpha\in (0,1]$:
\begin{equation}\label{Q-alpha-p}
Q_{\alpha,p}(v,u)=(p-1)\begin{pmatrix}
              \displaystyle{v^{\gamma}\left(1+\Gamma\left(\frac{u}{v}-1\right)\right)} & \displaystyle{-\frac{1}{2}(v^{\gamma}+u^{\gamma})}   \\
              \displaystyle{-\frac{1}{2}(v^\gamma+u^\gamma)} & \displaystyle{u^\gamma\left(1+\Gamma\left(\frac{v}{u}-1\right)\right)}
        \end{pmatrix} \quad \quad \mbox{with}\quad \Gamma = \frac{1-\alpha}{\alpha(p-1)}.
\end{equation}
and set $Q_{1,1}(v,u)=0$.
We can check that this matrix is nonnegative for all $u,v>0$ if and only if $(\alpha,p)\in K_{|n|}$ (see Lemma \ref{lem-sign-matrix-Q}).
For $y\in \R$, we set $y_+=\max(0,y)$. Then our main result is the following

\begin{theo}\label{theo-positive-negative-part}{\bf (The new contraction family)}\\
We recall that $\Omega=\R^d$ with $d\geq 1$. Let $|n|<1$ and $U$ (resp. $V$) be the solution of (\ref{eq::1}) in the sense of Definition \ref{defi::1}, with initial data $0\leq U_0\in L^1(\Omega)\cap L^\infty(\Omega)$ (resp. $0\leq V_0\in L^1(\Omega)\cap L^\infty(\Omega)$).
Let $(\alpha,p)\in K_{|n|}$ where $K_{|n|}$ is defined in (\ref{eq::2}). We set 
$$u=U^\alpha \quad \mbox{and} \quad v=V^\alpha.$$
Then for $w=u$ and $v$, we have $w\in C(Q)\cap L^\infty(Q)\cap C([0,+\infty);L^p(\Omega))$ and $w\in C^\infty(Q\cap \left\{w>0\right\})$.
Moreover, with the notation  $v(\tau)=v(\tau,\cdot)$, $u(\tau)=u(\tau,\cdot)$, we have for all $t\ge 0$
\begin{equation}\label{eq::4}
\int_{\Omega}\left(v(t)-u(t)\right)_{+}^p~dx + p \int_0^td\tau  \int_{\{v-u>0\}\cap \left\{u>0\right\}} e[v,u]\ dx \le \int_{\Omega}\left(v(0)-u(0)\right)_{+}^p~dx
\end{equation}
with
\begin{equation}\label{eq::r11}
e[v,u]=|v-u|^{p-2}
\begin{pmatrix}
\nabla v  \\
\nabla u
\end{pmatrix}^T
Q_{\alpha,p}(v,u)
\begin{pmatrix}
\nabla v \\
\nabla u
\end{pmatrix}\ge 0, \quad \mbox{for} \quad v>u>0.
\end{equation}
\end{theo}
Note that it is possible to see that there is equality in (\ref{eq::4}), for instance if $u$ and $v$ are smooth and positive on the torus  $\T^d$
instead of the whole space.

\noindent We now define a second symmetric matrix $M_{1,1}=0$ and for $(\alpha,p)\in K_{|n|}\backslash \left\{(1,1)\right\}$
\begin{equation}\label{K-alpha-p}
M_{\alpha,p}=(p-1)\begin{pmatrix}
              1 & \displaystyle{-\Gamma+\frac{\gamma}{2}}   \\
                  &      \\
              \displaystyle{-\Gamma+\frac{\gamma}{2}} & \displaystyle{\Gamma (1-\gamma)}
        \end{pmatrix},
\end{equation}
which turns out to be also nonnegative (see Lemma \ref{lem::r55}).\\

\noindent Then, we have also the following result
\begin{theo}\label{th::r3}{\bf (Decay of the gradient)}\\
Let $U$ be the solution of equation (\ref{eq::1}) with initial data $U_0$, under the assumptions of Theorem \ref{theo-positive-negative-part}.
Assume that $\nabla U_0^\alpha \in L^p(\Omega)$.  Setting $u=U^\alpha$, we have $\nabla u \in L^\infty([0,+\infty);L^p(\Omega))$.
Moreover, with the notation $u(\tau)=u(\tau,\cdot)$, we have for almost every $t\ge 0$:
\begin{equation}\label{eq::r50}
\begin{array}{l}
\displaystyle \int_{\Omega}\left|\nabla u(t)\right|^p~dx
+ p\int_0^t d\tau \int_{\{\nabla u\not= 0\}\cap \left\{u>0\right\}} {\bar e}[u, \nabla u,D^2 u]
dx\quad \le  \int_{\Omega}\left|\nabla u(0)\right|^p~dx.
\end{array}
\end{equation}
with
\begin{equation}\label{eq::r51}
\bar e[u, w, A]=
|w|^{p-2}u^{\gamma-2}
\begin{pmatrix}
u A_0 \\
w^2
\end{pmatrix}^T
M_{\alpha,p}
\begin{pmatrix}
uA_0\\
w^2
\end{pmatrix}+ |w|^{p-2} u^{\gamma}\left(|A_2|^2+p|A_1|^2\right)\ge 0,
\end{equation}
where   we write a general symmetric matrix $A$  as
$$A=\begin{pmatrix}
A_0 & A_1 \\
A_1^T & A_2
\end{pmatrix} \quad \mbox{by blocks in the vector space}\quad \R b \oplus b^\perp\quad \mbox{with}\quad b=\frac{w}{|w|}.$$
\end{theo}

\noindent Again, note that it is possible to see that there is equality in (\ref{eq::r50}), if $U_0$ is smooth and positive with $\Omega=\T^d$.\\

Let us mention that our results will be used in a future work to get information on the anomalous exponents associated to self-similar solutions like the one of Aronson-Graveleau (see \cite{AG}), and also the one of King-Peletier-Zhang (see \cite{K, PZ}).

\section{The new contraction property}\label{sec-proof-theo-wsol}

This section is devoted to the proof of Theorem \ref{theo-positive-negative-part}.
To this end, we first need the following result whose proof is elementary.
\begin{lem}{\bf(Properties of the matrix $Q_{\alpha,p}$)}\label{lem-sign-matrix-Q}\\
For $|n|<1$, let $\displaystyle (\alpha, p)\in \mbox{Int}\ K_{|n|}$ where $K_{|n|}$ is defined in (\ref{eq::2}).\\
\noindent {\bf i) Positive definite matrix}\\
Let $\eta >0$. Then for $v>u>0$ with $1+\eta^{-1}> v/u \ge 1+\eta$, we have
\begin{equation}\label{eq::r21}
Q_{\alpha,p}(v,u)\ge v^\gamma \nu I_2,
\end{equation}
where $I_2$ is the $2\times 2$ identity matrix, $Q_{\alpha,p}$ is defined in (\ref{Q-alpha-p}), and $\nu=\nu(n,\alpha,p,\eta)>0$.\\
\noindent {\bf ii) Bound from below for the top diagonal term}\\
Let $v\ge u>0$, then the top diagonal term of the matrix satisfies
\begin{equation}\label{eq::r27}
(Q_{\alpha,p}(v,u))_{11} \ge v^\gamma \nu_0,
\end{equation}
for some $\nu_0=\nu_0(n,\alpha,p)>0$.
\end{lem}

\noindent \textbf{Proof of Lemma \ref{lem-sign-matrix-Q}}\\
We set $\displaystyle{w=u/v \in (0,1]}$ and start with $\alpha\in (0,1)$ and $p>1$.\\
\textbf{Step 1: on the determinant}\\
A direct computation shows that $\mbox{det}\  Q_{\alpha,p}(v,u)> 0$ if and only if
\begin{equation}\label{eq::r25}
G_{\gamma}(w):=\Gamma\left(1-\Gamma\right)f_1(w)-\frac{1}{4}f_{\gamma}(w)> 0 \quad \mbox{with}\quad f_{\gamma}(w):=\displaystyle{w^{\gamma}}-2+\frac{1}{w^{\gamma}}.
\end{equation}
Note that $G_\gamma(1)=G'_\gamma(1)=0$. We now show that $G_\gamma$ is strictly convex for $w\le 1$, which will imply (\ref{eq::r25}) for $w\not=1$.
Indeed, we have
$$G_{\gamma}''(w)=-w^{-3}F_{\gamma}(w) \quad \mbox{with} \quad F_{\gamma}(w):=\frac{1}{4}\gamma\left((\gamma-1)w^{\gamma+1}+(\gamma+1)w^{1-\gamma}\right)-2\Gamma\left(1-\Gamma\right)\le  F_\gamma(1),$$
where we have used $|\gamma|< 1$ and $w\le 1$ to get the inequality. We see that 
$F_{\gamma}(1)< 0$ if and only if $\Gamma$ is in between the two distinct roots 
$\displaystyle{\Gamma_{\pm}:=\frac{1\pm \sqrt{1-\gamma^2}}{2}}$,
which is equivalent to $(\alpha, p)\in \mbox{Int}\ K_{|n|}$. This implies the strict convexity of $G_\gamma$.\\
\textbf{Step 2: conclusion}\\
Recall that for $(\alpha,p)\in \mbox{Int}\ K_{|n|}$, we have $\displaystyle p> P_-(\alpha)\ge \frac{1}{\alpha}$.
We deduce that each diagonal term of $Q_{\alpha, p}$ is positive (which shows ii)) and then $\mbox{tr}\ Q_{\alpha, p}(v,u)> 0$.
With Step 1, this implies the result i).
\hfill $\Box$\\

\noindent {\bf Proof of Theorem \ref{theo-positive-negative-part}}\\
We only deal with the case $(\alpha,p)\in \mbox{Int}\ K_{|n|}$. Note that the border case can easily be recovered,
by a passage to the limit.
The proof is divided into several steps.\\
\noindent {\textbf{Step 1: Estimate for smooth positive solutions}}\\
Given two nonnegative and smooth initial data $U_0, V_0$ with compact support,
we consider modified initial data for $\varepsilon\in (0,1)$
$$U_{0\varepsilon}= U_0+\varepsilon \quad \mbox{and} \quad  V_{0\varepsilon}=V_0+\varepsilon.$$
According to the standard parabolic theory (see \cite{LSU}), 
we can consider the smooth positive functions  $U_\varepsilon$, $V_\varepsilon$, which are the classical solutions of (\ref{eq::1}) respectively associated to initial data $U_{0\varepsilon}$ and $V_{0\varepsilon}$. We have moreover the obvious bounds
\begin{equation}\label{eq::r22}
\varepsilon\le U_\varepsilon, V_\varepsilon \le |U_0|_{L^\infty(\Omega)} + \varepsilon.
\end{equation}
For $u_\varepsilon=U_\varepsilon^\alpha$, $v_\varepsilon = V_\varepsilon^\alpha$, for a cut-off function $\varphi\in C_c^\infty(\Omega)$ satisfying $\varphi\ge 0$ and for $p>1$, we will compute $\frac{d}{dt}\int_\Omega \varphi^2 (v_\varepsilon-u_\varepsilon)^p_+$.
We will use the positivity properties of the matrix $Q_{\alpha,p}$, to control the error term created by the cut-off. To this end,  it is useful to introduce for $\delta \in [0,\varepsilon)$, the function $\Psi_{\delta,p}(w)$ which is an approximation of $\Psi_{0,p}(w)=w^p_+$, given by
$$\Psi_{\delta,p}(w)=\left\{\frac{w^p}{p}-\delta^{p-1}w-\delta^p\left(\frac1p -1\right)\right\}\cdot 1_{\left\{w>\delta\right\}}.$$
Then a direct computation gives
\begin{equation}\label{eq::r12}
\begin{array}{l}
\displaystyle \frac{d}{dt}\int_{\Omega} \varphi^2 \Psi_{\delta,p}\left(v_\varepsilon-u_\varepsilon\right)dx\\
\\
\displaystyle =   \int_{\Omega} \varphi^2 \Psi_{\delta,p}'\left(v_\varepsilon-u_\varepsilon\right)
\left\{v_\varepsilon^{\gamma}\Delta v_\varepsilon-u_\varepsilon^{\gamma}\Delta u_\varepsilon\right\} \ dx
+ \int_{\Omega}\varphi^2 \Psi_{\delta,p}'\left(v_\varepsilon-u_\varepsilon\right)\bar \gamma  
\left\{v_\varepsilon^{\gamma-1}\vert \nabla v_\varepsilon\vert^2  -u_\varepsilon^{\gamma-1}\vert \nabla u_\varepsilon\vert^2 \right\}\ dx\\
\displaystyle =
 \int_{\{v_\varepsilon-u_\varepsilon>\delta\}} f\ dx
-\int_{\{v_\varepsilon-u_\varepsilon>\delta\}}
\varphi^2 e[v_\varepsilon,u_\varepsilon]\ dx
+ \delta^{p-1}\int_{\left\{v_\varepsilon-u_\varepsilon>\delta\right\}} \varphi^2 g \ dx
\end{array}
\end{equation}
with the notation $e[\cdot,\cdot]$ defined in (\ref{eq::r11}) and
$$f= -\Psi_{\delta,p}'\left(v_\varepsilon-u_\varepsilon\right)2\varphi \nabla \varphi \cdot
\left\{v_\varepsilon^{\gamma}\nabla v_\varepsilon-u_\varepsilon^{\gamma}\nabla u_\varepsilon\right\}  
\quad \mbox{and}\quad g = \frac{(\alpha-1)}{\alpha}\left(|\nabla v_\varepsilon|^2 v_\varepsilon^{\gamma-1} - |\nabla u_\varepsilon|^2 u_\varepsilon^{\gamma-1}\right).$$
Here in (\ref{eq::r12}) we have used equation (\ref{eq::1-bis}) to get the second line, and have done an integration by parts on the first term of the second line, using in particular for $\delta>0$ the chain rule 
$\nabla \Psi_{\delta,p}'\left(w\right) = \Psi_{\delta,p}''\left(w\right) \nabla w$
with equality almost everywhere with the convention that the right hand side is zero if $\nabla w=0$ irrespective of whether $\Psi_{\delta,p}''$ is defined. The case $\delta=0$ can be recovered, passing to the limit $\delta\to 0$.\\
Using bounds (\ref{eq::r22}), we see that $0<v_\varepsilon^\gamma, u_\varepsilon^\gamma \le M$ for some  constant $M>0$. Moreover for $v_\varepsilon-u_\varepsilon>\delta>0$, we have by Young inequality
\begin{equation}\label{eq::r33}
|f| \le (v_\varepsilon-u_\varepsilon)^{p-2}\left\{\frac{2M}{\varepsilon\nu} (v_\varepsilon-u_\varepsilon)^2 |\nabla \varphi|^2 + \varphi^2 \varepsilon\nu (v_\varepsilon^\gamma |\nabla v_\varepsilon|^2 +u_\varepsilon^\gamma |\nabla u_\varepsilon|^2)\right\}
\end{equation} 
where $\nu=\nu(\delta)>0$ is given in (\ref{eq::r21}) for $\delta>0$.\\
\noindent {\textbf{Step 2: First integral estimate on the gradient}}\\
We now choose $u_\varepsilon\equiv \varepsilon$ and $\delta=0$.
Then estimate (\ref{eq::r33}) still holds true, but with $\nu$ replaced by $\nu_0$  given in (\ref{eq::r27}). This implies
$$\displaystyle \frac{d}{dt}\int_{\Omega} \varphi^2 \Psi_{0,p}\left(v_\varepsilon-\varepsilon\right)dx
\le  \int_{\Omega} C_0 |\nabla \varphi|^2 \Psi_{0,p}(v_\varepsilon-\varepsilon)\ dx
-(1-\varepsilon)\int_{\{v_\varepsilon-\varepsilon>0\}}
\varphi^2 e[v_\varepsilon,\varepsilon]\ dx$$
with $\displaystyle C_0=\frac{2pM}{\varepsilon\nu_0}$.
We then take a sequence of functions $\varphi$ converging towards $\varphi_\lambda(x)=e^{-\lambda |x|}$ for $\lambda>0$. Integrating in time, and passing to the limit $\lambda\to 0$, and using again (\ref{eq::r27}), we get
\begin{equation}\label{eq::r28}
\int_{\Omega} \Psi_{0,p}\left(v_\varepsilon(t)-\varepsilon\right)dx+ (1-\varepsilon)\nu_0 \int_0^t d\tau \int_{\Omega}  v_\varepsilon^\gamma |\nabla v_\varepsilon|^2 \le \int_{\Omega} \Psi_{0,p}\left(v_\varepsilon(0)-\varepsilon\right)dx
\end{equation}

\noindent {\textbf{Step 3:  Refined estimate}}\\
We come back to general $v_\varepsilon,u_\varepsilon$ as in Step 1 and now consider $\delta>0$.
Using (\ref{eq::r33}), we get
$$\begin{array}{l}
\displaystyle \frac{d}{dt}\int_{\Omega} \varphi^2 \Psi_{\delta,p}\left(v_\varepsilon-u_\varepsilon\right)dx\\
\\
\displaystyle\le  \int_{\Omega} C |\nabla \varphi|^2 \Psi_{0,p}(v_\varepsilon-u_\varepsilon)\ dx
-(1-\varepsilon)\int_{\{v_\varepsilon-u_\varepsilon>\delta\}}
\varphi^2 e[v_\varepsilon,u_\varepsilon]\ dx
+\delta^{p-1}\int_{\{v_\varepsilon-u_\varepsilon>\delta\}} \varphi^2 g\ dx
\end{array}$$
with $\displaystyle C=\frac{2pM}{\varepsilon\nu}$.
Note that estimate (\ref{eq::r28})  for both $v_\varepsilon$ and $u_\varepsilon$, controls uniformly in time $\int_\Omega \Psi_{0,p}(v_\varepsilon-u_\varepsilon)\ dx$ and the time integral of $\int_\Omega \delta^{p-1}|g|\ dx$.
Therefore, we can again apply the choice of $\varphi=\varphi_\lambda$, integrate in time, and take the limit $\lambda\to +\infty$ as in Step 2, and then conclude  in the limit $\delta\to 0$
\begin{equation}\label{eq::r30}
\displaystyle \int_{\Omega} \Psi_{0,p}\left(v_\varepsilon(t)-u_\varepsilon(t)\right)dx
+ (1-\varepsilon) \int_0^td\tau\ \int_{\{v_\varepsilon-u_\varepsilon>0\}\cap \left\{u_\varepsilon>0\right\}}
e[v_\varepsilon,u_\varepsilon]\ dx \le \int_{\Omega} \Psi_{0,p}\left(v_\varepsilon(0)-u_\varepsilon(0)\right)dx.
\end{equation}

\noindent {\textbf{Step 4: The limit $\varepsilon\to 0$}}\\
For a point $P_0=(t_0,x_0)$,  let us denote the open parabolic cylinder $Q_r(P_0)=(t_0-r^2,t_0)\times B_r(x_0)$,
where $B_r(x_0)$ is the open ball of center $x_0$ and radius $r>0$.
It is known by Theorem 1.1 in Sacks \cite{S83} that for any smooth solution $U_\varepsilon$ of (\ref{eq::1}) on $Q_{2r}(P_0)$, there exists a modulus of continuity $\omega$ of $U_\varepsilon$ on $Q_{r}(P_0)$, depending only on $r,d,n$ and $|U_\varepsilon|_{L^\infty(Q_{2r}(P_0))}$.
This property is automatically transfered to $u_\varepsilon=U_\varepsilon^\alpha$
with the modulus of continuity $\omega^\alpha$. 
This implies, by Ascoli-Arzela theorem, that $U_\varepsilon\to U$, where $U$ is still a solution of (\ref{eq::1}) with initial data $U_0$
in the sense of Definition \ref{defi::1}. We also note that (with $\alpha p\ge 1$)
$$\int_\Omega |u(t,\cdot)-u(s, \cdot)|^p \le \int_\Omega |U(t,\cdot)-U(s, \cdot)|^{\alpha p} \le (2|U|_{L^\infty(Q)})^{\alpha p -1} |U(t,\cdot)-U(s, \cdot)|_{L^1(\Omega)}$$
which shows in particular that $U\in C([0,+\infty);L^1(\Omega))$ implies $u\in C([0,+\infty);L^p(\Omega))$.

From the standard parabolic theory  \cite{LSU}, we know that  $U\in C^\infty(Q\cap \left\{U>0\right\})$, with corresponding quantitative estimates on solutions locally bounded from above and below.
Similarly, we have $V_\varepsilon\to V$, and we call $u=U^\alpha$, $v=V^\alpha$.
Using the $C^1$ convergence of $(v_\varepsilon,u_\varepsilon)$ to $(v,u)$ on compact sets inside $\left\{v-u >0\right\}\cap \left\{u> 0\right\}$, we can pass to the limit in (\ref{eq::r30}), and get the same inequality for $\varepsilon=0$ with $(v_\varepsilon,u_\varepsilon)$ replaced by $(v,u)$.
This means (\ref{eq::4}).
We finally conclude to the result for general initial data $U_0,V_0$, by a standard approximation argument in $L^1(\Omega)$.
\hfill $\Box$

\section{Decay of the gradient}\label{sec-gradient}

This section is divided into two subsections. In the first subsection, we give the proof of the gradient decay Theorem \ref{th::r3}.
In the second subsection, we give a directional derivative estimate (Theorem \ref{theo-d-derivative}) as a corollary of our contraction estimate.

\subsection{Decay of the gradient}\label{s4.1}

We start with the following simpler analogue of Lemma \ref{lem-sign-matrix-Q}, whose proof follows from an elementary computation.

\begin{lem}{\bf(Properties of the matrix $M_{\alpha,p}$)}\label{lem::r55}\\
For $|n|<1$, let $\displaystyle (\alpha, p)\in \mbox{Int}\ K_{|n|}$ where $K_{|n|}$ is defined in (\ref{eq::2}).\\
Then we have 
\begin{equation}\label{eq::r56}
M_{\alpha,p}\ge \nu_1 I_2,
\end{equation}
where $I_2$ is the $2\times 2$ identity matrix, $M_{\alpha,p}$ is defined in (\ref{K-alpha-p}), and $\nu_1=\nu_1(n,\alpha,p)>0$.
\end{lem}

\noindent {\bf Proof of Theorem \ref{th::r3}}\\
We follow the lines of proof of Theorem \ref{theo-positive-negative-part}.
For $\delta>0$, a direct computation gives with $\displaystyle b_\varepsilon=\frac{\nabla u_\varepsilon}{|\nabla u_\varepsilon|}$
\begin{equation}\label{eq::r53}
\begin{array}{l}
\displaystyle \frac{d}{dt}\int_{\Omega} \varphi^2 \Psi_{\delta,p}\left(|\nabla u_\varepsilon|\right)dx\\
\\
\displaystyle =   \int_{\Omega} \varphi^2 \Psi_{\delta,p}'\left(|\nabla u_\varepsilon|\right)
b_\varepsilon \cdot \nabla \left\{u_\varepsilon^{\gamma-1}\left(u_\varepsilon \Delta u_\varepsilon +\bar \gamma |\nabla u_\varepsilon|^2 \right)\right\} \ dx\\
\displaystyle =
 \int_{\{|\nabla u_\varepsilon|>\delta\}} {\bar f}\ dx
-\int_{\{|\nabla u_\varepsilon|>\delta\}}
\varphi^2 {\bar e}[u_\varepsilon,\nabla u_\varepsilon, D^2 u_\varepsilon]\ dx
+\delta^{p-1} \int_{\{|\nabla u_\varepsilon|>\delta\}}
\varphi^2  {\bar g}
\end{array}
\end{equation}
with the notation $\bar{e}[\cdot,\cdot,\cdot]$ defined in (\ref{eq::r51}) and
$${\bar f}  = \sum_{i=1}^d -2\varphi \nabla_i \varphi \Psi_{\delta,p}'\left(|\nabla u_\varepsilon|\right)
b_\varepsilon \cdot \nabla \left\{u_\varepsilon^{\gamma} \nabla_i u_\varepsilon \right\}
\quad \mbox{and}\quad {\bar g} = \frac{(\alpha-1)}{\alpha} b_\varepsilon \cdot \nabla \left\{u_\varepsilon^{\gamma-1}|\nabla u_\varepsilon|^2\right\}.$$
Here we have integrated by parts the term $\Delta u_\varepsilon$ in the second line of (\ref{eq::r53}).
We then pass to the limit $\delta\to 0$. We replace the fine property (\ref{eq::r21}) of the matrix $Q_{\alpha,p}$, by (\ref{eq::r56}). We finally conclude as in the proof of Theorem \ref{theo-positive-negative-part}. \hfill $\Box$

\subsection{Decay of the directional derivative}\label{s4.2}

For $\xi=(\xi^t,\xi^x)\in \R\times \R^d$, we define the directional derivative
\begin{equation}\label{eq::5}
D_\xi u=\xi^t u_t + \xi ^x\cdot \nabla u.
\end{equation}

In this subsection, we prove the following result as a consequence of our contraction estimate.\\

\begin{theo}\label{theo-d-derivative}{\bf (Decay of the directional derivative)}\\
Let $U$ be the solution of equation (\ref{eq::1}) with initial data $U_0$, under the assumptions of Theorem \ref{theo-positive-negative-part}.
We also set $u=U^\alpha$. 
Let $\xi \in  \R\times \R^d$, and let us assume that
\begin{equation}\label{eq::r10}
C_\xi := \limsup_{\eta\to 0^+} \int_{\Omega} \left|\frac{u((0,x)+ \eta\xi)-u(0,x)}{\eta}\right|^p\ dx < +\infty
\end{equation}
With the notation $D_\xi u$ defined in (\ref{eq::5}), we have $D_\xi u \in L^\infty([0,+\infty);L^p(\Omega))$.
Moreover, with the notation $u(\tau)=u(\tau,\cdot)$, we have for almost every $t\ge 0$:
\begin{equation}\label{contraction-deriv}
\begin{array}{l}
\displaystyle \int_{\Omega}\left|D_{\xi}u(t)\right|^p~dx
+ p\int_0^t\int_{\{D_{\xi}u\not= 0\}\cap \left\{u>0\right\}} \bar {\bar e}[u,D_\xi u,\nabla u]
dx\quad \le C_\xi.
\end{array}
\end{equation}
with
\begin{equation}\label{eq::r20}
\bar {\bar e}[u,w,\nabla u]=
|w|^{p-2}u^{\gamma-2}
\begin{pmatrix}
u\nabla w  \\
w\nabla u
\end{pmatrix}^T
M_{\alpha,p}
\begin{pmatrix}
u\nabla w  \\
w\nabla u
\end{pmatrix}\ge 0.
\end{equation}

\end{theo}
In particular (\ref{eq::r10}) is satisfied if $\xi^t=0$ and $D_\xi U_0^\alpha \in L^p(\Omega)$.

\noindent We will use the following result.
\begin{lem}\label{lem-lim-Q-alpha-p}{\bf (Limit of the quadratic term)}\\
For $|n|<1$, let $\displaystyle (\alpha, p)\in K_{|n|}$ where $K_{|n|}$ is defined in (\ref{eq::2}).\\
Let $\Omega_1$ be an open set such that $u,w,z\in C^1(\Omega_1)$ with $u>0$ on $\Omega_1$.
For $\eta \in (0,1)$, let $v_\eta$ be a function satisfying
$$|v_\eta -(u+\eta w + \eta^2 z)|_{C^1(\Omega_1)} =o(\eta^2)$$
Then we have on $\Omega_1$
$$J_\eta(v_\eta, u) 
= \eta^{-2 }\begin{pmatrix}
\nabla v_\eta \\
\nabla u
\end{pmatrix}^T
Q_{\alpha,p}\left(v_\eta,u\right)
\begin{pmatrix}
\nabla v_\eta\\
\nabla u
\end{pmatrix}
\quad \underset{\eta \to 0}{\longrightarrow}\quad  J_0=u^{\gamma-2}
\begin{pmatrix}
u\nabla w  \\
w\nabla u
\end{pmatrix}^T
M_{\alpha,p}
\begin{pmatrix}
u\nabla w \\
w\nabla u
\end{pmatrix},$$
where we recall that the matrices $Q_{\alpha, p}$ and $M_{\alpha, p}$ are respectively defined  in (\ref{Q-alpha-p}) and (\ref{K-alpha-p}).
\end{lem}

\noindent \textbf{Proof of Lemma \ref{lem-lim-Q-alpha-p}}\\
The proof is done using a simple Taylor expansion argument.
For $Q_{\alpha,p}(v,u)$, we set $Q_i = \partial^i_v Q_{\alpha,p}(u,u)$ for $i=0,1,2$.
For $V_\eta=\begin{pmatrix}
\nabla v_\eta  \\
\nabla u
\end{pmatrix}$, we set $V_\eta = V_0+ \eta V_1 +\eta^2 V_2 + o(\eta^2)$.
Using the fact that $Q_0 V_0=0 = V_0^TQ_1V_0$, we get
$$J_\eta(v_\eta,u) = \frac12 w^2 V_0^T Q_2  V_0 + 2w V_1^T Q_1 V_0 + V_1^T Q_0 V_1 + o(1)$$
which implies the result.
\hfill $\Box$\\

\noindent {\bf Proof of Theorem \ref{theo-d-derivative}}\\
Let $v_\eta(t,x)=u((t,x)+\eta \xi)$.
We simply apply Theorem \ref{theo-positive-negative-part} to $v_\eta-u$ and $u-v_\eta$ and use the symmetry $e[v,u]=e[u,v]$,
in order to control
$$\int_\Omega \left|\frac{v_\eta(t)-u(t)}{\eta}\right|^p$$
We pass to the limit as $\eta\to 0$ using Lemma \ref{lem-lim-Q-alpha-p} which gives (\ref{eq::r10}).
\hfill $\Box$

\noindent {\bf Aknowledgements}


\noindent The authors thank B. Andreianov for helpful indications on the literature.



\end{document}